\documentclass{article}

\usepackage{amsmath,amssymb,enumerate,bbm,mathrsfs}

\newtheorem{cntr}{ERROR! SHOULD NOT USE THIS}

\newcommand{\tmop}[1]{\operatorname{#1}}

\newtheorem{varremark}[cntr]{Remark}
\newenvironment{remark}{\begin{varremark}\em}{\em\end{varremark}}
\newcommand{\nin}{\not\in}

\newtheorem{varnote}{Note}

\newenvironment{proof*}[1]{
  \noindent\textbf{#1\ }}{\hspace*{\fill}
  \begin{math}\Box\end{math}\medskip}

\newtheorem{theorem}[cntr]{Theorem}

\newtheorem{algo}{Algorithm}

\newcommand{\constfont}[1]{\mathbf{#1}}

\newcommand{\psia}[0]{\Psi}

\newcommand{\xs}[0]{x_0}
\newcommand{\ks}[0]{k_0}
\newcommand{\sdev}[0]{\sigma}
\newcommand{\oversample}[0]{\constfont{M}}

\newcommand{\frmlt}[3]{\phi_{(\vec{#1},\vec{#2})}(\vec{#3})}

\newcommand{\frmltDual}[3]{\tilde{\phi}_{(\vec{#1},\vec{#2})}(\vec{#3})}

\newcommand{\frmltGen}[2]{\phi_{#1}(#2)}

\newcommand{\frmltDef}[3]{ \pi^{-N/4}\sdev^{-N/2} e^{i \ks \vec{#2} \cdot \vec{#3}} e^{- \absSmallTwo{ \vec{#3} - \vec{#1} \xs }^2 / 2 \sdev^2}}

\newcommand{\dualW}[0]{\tilde{g}}
\newcommand{\frmltDualDef}[3]{e^{i \ks \vec{#2} \cdot \vec{#3}}\dualW(\vec{#3} -\vec{ #1} \xs)}

\newcommand{\frmltProp}[4]{\frac{
    \exp \left(i \vec{#2} \ks \cdot(\vec{#3} - (\vec{#2}\ks/2) #4 - \vec{#1} \xs) \right)
  }{\pi^{N / 4} \sdev^{N/2} ( 1 + i #4 / \sdev^2 )^{N / 2}}
  \exp \left( \frac{- \absSmallTwo{ \vec{#3} - \vec{#2} \ks t - \vec{#1} \xs}^2}{2 \sdev^2 ( 1 + i #4 / \sdev^2 )} \right)
}

\newcommand{\fcoeffs}[3]{ {#1}_{(\vec{#2},\vec{#3})}}

\newcommand{\ZtwoN}[0]{\mathbb{Z}^N \times \mathbb{Z}^N}
\newcommand{\Zn}[0]{{\mathbb{Z}^N}}
\newcommand{\Rn}[0]{{\mathbb{R}^N}}
\newcommand{\fsum}[3]{\sum_{(\vec{#1},\vec{#2}) \in #3}}
\newcommand{\fsumGen}[2]{\sum_{#1 \in  #2}}

\newcommand{\proj}[4]{P^{#3}_{#1; #2}(#4)}

\newcommand{\psl}[1]{\mathcal{P}_{#1}}

\newcommand{\frmLB}[0]{A_F}
\newcommand{\frmUB}[0]{B_F}

\newcommand{\DWdecayConstXd}[1]{\constfont{g}(\xs,\ks,N,#1)}

\newcommand{\DWdecayRateX}[0]{\constfont{r}(\xs,\ks)}

\newcommand{\XBuffCube}[3]{\constfont{X}}
\newcommand{\KBuffCube}[2]{\constfont{K}}

\newcommand{\Gammai}{\Gamma^{-1}}

\newcommand{\frmltbbDef}[3]{ UNDEFINED AS OF YET}

\newcommand{\kmax}[0]{k_{\tmop{max}}}
\newcommand{\kmin}[0]{k_{\tmop{min}}}
\newcommand{\Tstep}[0]{T_{\tmop{step}}}
\newcommand{\Tmax}[0]{T_{\tmop{max}}}

\newcommand{\kmaxNL}[0]{k_{\tmop{max},\tmop{NL}}}
\newcommand{\LNL}[0]{L_{\tmop{NL}}} 

\newcommand{\Lfree}[0]{L_{F}}

\newcommand{\vx}[0]{\vec{x}}

\newcommand{\vk}[0]{\vec{k}}
\newcommand{\va}[0]{\vec{a}}
\newcommand{\vb}[0]{\vec{b}}

\newcommand{\Lb}[0]{L_{\tmop{int}}}
\newcommand{\wb}[0]{w}

\newcommand{\LbC}[0]{L_{\tmop{comp}}}

\newcommand{\IBox}[0]{[-\Lb,\Lb]^N}
\newcommand{\FBox}[0]{[-(\Lb+\wb),(\Lb+\wb)]^N}
\newcommand{\CBox}[0]{[-\LbC,\LbC]^N}

\newcommand{\LNLBox}[0]{[-\LNL,\LNL]^N}

\newcommand{\norm}[2]{\left\| #1 \right \|_{#2}}
\newcommand{\abs}[1]{\left| #1 \right|}

\newcommand{\absSmallN}[2]{| #1 |_{#2}}

\newcommand{\absSmallOne}[1]{\absSmallN{#1}{1}}
\newcommand{\absSmallTwo}[1]{\absSmallN{#1}{2}}

\newcommand{\absSmallInfty}[1]{\absSmallN{#1}{\infty}}

\newcommand{\Hs}[0]{H^s}
\newcommand{\Hsb}[0]{H^s_b}
\newcommand{\IP}[2]{\left\langle #1 | #2 \right\rangle}

\newcommand{\nlin}[2]{g(#1,\vec{x},#2) #2 }
\newcommand{\nlinsemi}[2]{g(#1,\vec{x},#2)  }

\newcommand{\nlinNoArg}[1]{g(#1,\vec{x},\cdot)}

\newcommand{\Lap}[0]{(1/2)\Delta}

\newcommand{\U}[1]{\mathcal{U}(#1)}

\newcommand{\Uf}[1]{e^{i \Lap #1}}

\numberwithin{algo}{section}
\numberwithin{equation}{section}
\numberwithin{cntr}{section}

\newcommand{\comment}[1]{}

\begin{document}

\title{Open Boundaries for the Nonlinear Schr\"odinger Equation}

\author{A. Soffer and C. Stucchio}

\maketitle

\begin{abstract}
  We present a new algorithm, the Time Dependent Phase Space Filter (TDPSF) which is used to solve time dependent Nonlinear Schrodinger Equations (NLS). The algorithm consists of solving the NLS on a box with periodic boundary conditions (by any algorithm). Periodically in time we decompose the solution into a family of coherent states. Coherent states which are outgoing are deleted, while those which are not are kept, reducing the problem of reflected (wrapped) waves. Numerical results are given, and rigorous error estimates are described.

  The TDPSF is compatible with spectral methods for solving the interior problem. The TDPSF also fails gracefully, in the sense that the algorithm notifies the user when the result is incorrect. We are aware of no other method with this capability. 
\end{abstract}

\section{Introduction and Definitions}

Consider a semilinear Schr\"odinger equation on $\mathbb{R}^{N+1}$
\begin{equation}
  \label{eq:NLSE}
  i \partial_t \psi ( x, t ) = - \Lap \psi ( x, t ) + \nlin{t}{\psi(\vx,t)}
\end{equation}
where $\nlinNoArg{t}$ is some semilinear perturbation. For instance, $\nlinNoArg{t}$ could be $V(\vx,t) + f(\abs{\psi(\vx,t)}^2)$ for some smooth function $f$ and (spatially) localized potential $V(\vx,t)$. Abusing terminology, we refer to $\nlinNoArg{t}$ as the nonlinearity even when it is linear. 

We assume the initial condition and nonlinearity are such that the nonlinearity remains localized inside some box $\LNLBox$. Outside this region the solution is assumed to behave like a free wave.

One very common method of solving such a problem is domain truncation. That is, one solves the PDE \eqref{eq:NLSE} numerically on a region $[-L,L]^N$. On the finite domain, boundary conditions must be specified. Dirichlet or Neumann boundaries introduce spurious reflections, while periodic boundaries allow outgoing waves to wrap around the computational domain. This causes the numerical solution to become incorrect after a time $T \approx L/\kmax$, where $\kmax$ is the ``maximal velocity'' (the highest relevant spatial frequency) of the solution.

There are two major approaches to dealing with this problem. The Dirichlet-to-Neumann method consists of attempting to use the exact solution as a boundary condition. This works rather well for the wave equation \cite{MR596431,MR658635,MR0471386,MR517938,776449}, but it is fraught with problems for dispersive waves and only limited progress has been made \cite{MR1924419,MR1785966,MR1869342,MR2114289}. Dirichlet-to-Neumann methods also prevent the use of fast spectral methods for solving the interior problem because spectral methods based on the FFT naturally impose periodic boundary conditions. A spectral solver for the interior problem is highly desirable since dispersive waves are difficult to calculate by other methods (see Remark \ref{remark:commentsOnPML}).

The other approach is to add a dissipative term which is localized on a buffer region $\FBox \setminus \IBox$. The dissipative term can be an absorbing potential \cite{neuhauser:complexPotentials} or the more sophisticated PML \cite{MR1294924}. This dissipates outgoing waves, but also dissipates incoming waves located near the boundary. This can introduce spurious dissipation (see Section \ref{sec:testsMediumRange}) which does not occur for the TDPSF.

\subsection{Our Approach}

We propose a new method. Suppose we want to approximate the solution in the space $H^{s}$ (with $s$ large enough to make \eqref{eq:NLSE} is well posed). We make the assumption that outside some region $[-\Lfree,\Lfree]^N$, $\psi(x,t)$ behaves like a free wave, i.e. $\psi(x,t) \approx e^{i \Lap t}\psi_+(x)$. The size of the computational box $\Lb$ must be larger than $\Lfree$.

We assume the nonlinearity is nearly localized both in position and momentum. That is, we assume the existence of $\LNL$, $\kmaxNL$ so that:
\begin{subequations}
\label{eq:defOfLnlKmaxNL}
  \begin{eqnarray}
    \norm{ [1-\chi_{[-\LNL,\LNL]^N}(\vx)] \nlin{t}{\psi(x,t)}}{\Hs} & \approx & 0 \\
    \norm{[1- \chi_{[-\kmaxNL,\kmaxNL]^N}(\vk)]\nlin{t}{\psi(x,t)}}{\Hs} & \approx & 0
  \end{eqnarray}
\end{subequations}
For the rigorous proof, we assume that the nonlinearity is Lipschitz in $H^{s}$, though this is probably unnecessary. In particular, the nonlinearity could take the form of a time dependent short range potential $V(\vx,t)\psi(\vx,t)$ or a local nonlinearity $f(\abs{\psi(\vx,t)}^{2 \sigma}) \psi(\vx,t)$ with $f(z)$ a bounded function.

We also assume that the solution remains localized in frequency, that is the $\Hs$ norm of  $\hat{\psi}(\vk,t)$ is small outside $[-\kmax,\kmax]^N$ for some large number $\kmax$ (the maximal momentum of the problem, which we assume exists).

Our algorithm is as follows. We assume the initial data is localized on a region $[-\Lb,\Lb]^N$. We solve \eqref{eq:NLSE} on the box $[-(\Lb+\wb),\Lb+\wb]^N$ on the time interval $[0,\Tstep]$. We insist that $\Lb > \Lfree, \LNL$ so that the solution behaves like a free wave on  $\FBox \setminus \IBox$.

By making $\Tstep$ smaller than $\wb/3\kmax$, we can ensure that $\psi(\vx,\Tstep)$ is mostly localized inside box $[-(\Lb+\wb),(\Lb+\wb)]^N$. Since very little mass has reached the boundary, the error is nearly zero \cite[Thm. 4.1]{us:TDPSFrigorous}.

We then decompose $\psi(\vx,\Tstep)$ into a sum of Gaussians (indexed by $\va,\vb \in \mathbb{Z}^N$, with lattice spacing $\xs$ in position, $\ks$ in momentum and standard deviation $\sdev$):
\begin{equation*}
  \psi(x,\Tstep)\approx \sum_{
    \begin{subarray}{l}
      \absSmallInfty{\va \xs} \leq \Lb+\wb\\
      \absSmallInfty{\vb \ks} \leq \kmax
    \end{subarray}
  } \fcoeffs{\psi}{a}{b}\frmltDef{a}{b}{x}
\end{equation*}
We then examine the Gaussians near the boundary, $\absSmallInfty{\va \xs} \geq \Lb$, and determine whether they are leaving the box or not under the free flow ($\absSmallN{\vx}{p}$ denotes the $l^p$ norm on finite vectors).  By assumption, $\Uf{t}$ is a sufficiently accurate approximation to the true solution for this part of the solution. We apply this approximation to each framelet:
\begin{multline}
  \label{eq:freeFrameletProp}
  \Uf{t} \frmltDef{a}{b}{x}\\
  = \frmltProp{a}{b}{x}{t}
\end{multline}
Under the free flow, a framelet moves along the trajectory $\va \xs + \vb \ks t$ while spreading about it's center. If a given Gaussian is leaving the box, we delete it (set $\fcoeffs{\psi}{a}{b}=0$). Some Gaussians spread more quickly than their center of mass moves, and we do not present here an algorithm to deal with these Gaussians. 

After this filtering operation, the only Gaussians remaining are either inside the box $\IBox$, or inside the box $\FBox$ but moving towards $\IBox$. Thus, it is safe to propagate for a time $\Tstep$, since what remains will not hit the boundaries before this time. We do the same at time $2\Tstep, 3\Tstep,$ etc.

The main drawback to this approach is that some Gaussians are ambiguous. Consider a Gaussian with velocity $0$: $\pi^{-N/4} \sdev^{-N/2} e^{-\absSmallTwo{\vx - \va \xs}^2/2\sdev^2}$. Suppose also that $\va \xs$ is located inside $\FBox \setminus \IBox$. By examining \eqref{eq:freeFrameletProp}, we observe that this framelet is spreading out laterally both into the box and outward. If we delete it, we have removed waves which should have returned. If we fail to delete it, then that part of the framelet which is spreading outwards will wrap around and cause an error.

This happens only for Gaussians which are moving slowly relative to the box. Thus, we impose one additional assumption -- we assume that Gaussians moving so slowly that they move both ways do not occur in the solution. The algorithm to deal with low frequency waves is more involved and we relegate this to a future work \cite{us:multiscale}.

\subsection{Error Bounds}

We prove rigorous error bounds for the TDPSF algorithm in \cite{us:TDPSFrigorous}, as well as stating explicitly the assumptions and defining what we mean by ``$\approx$''. Although the error bound is too long to state here, we describe briefly the form it takes.

For a general time-stepping algorithm (with periodic boundaries and no filtering), the error bound would take the following form:
\begin{multline}
  \label{eq:simpleErrorFormula}
  \sup_{t \in [0,\Tmax]} \norm{\U{t} \psi_0(x) - \psia(x,t)}{\Hsb} \leq
  \operatorname{BoundaryError}(\Tmax)\\
  + \operatorname{High Frequency}(\Tmax) + \operatorname{Low Frequency}(\Tmax)\\
  + \operatorname{Nonlocal Nonlinearity}(\Tmax) + \operatorname{Instability}(\Tmax)
\end{multline}

The term $\operatorname{BoundaryError}(\Tmax)$ encompasses errors due to waves wrapping/reflecting from the boundaries of the box. For many problems, this is the dominant error term. It is directly proportional to the mass which would have (if we were solving the problem on $\Rn$) radiated outside the box $\IBox$.

The $\operatorname{HighFrequency}(\Tmax)$ part stems from waves with momenta too high to be resolved by the discretization. The term $\operatorname{LowFrequency}(\Tmax)$ encompasses errors due to waves with wavelength that is long in comparison to the box. The term $\operatorname{Nonlocal Nonlinearity}(\Tmax)$ stems from that fraction of the nonlinearity itself which is located outside the box. The $\operatorname{Instability}(\Tmax)$ stems from the possibility that the dynamics of the solution itself might amplify the other errors dramatically (e.g. in strongly nonlinear problems).

Our algorithm reduces the term $\operatorname{BoundaryError}(\Tmax)$. The other errors are assumed to be small, since we are concerned only with the boundary error. 

The method of proof is a direct calculation. We calculate the errors made in the filtering step as well as the errors made while propagating in between filtering. We then add up the error over time, taking into account instabilities of the system being simulated. The error takes the form described in \eqref{eq:simpleErrorFormula}. The $\operatorname{BoundaryError}(\Tmax)$ term can be reduced by increasing $\wb$, and the error behaves like $\operatorname{BoundaryError}(\Tmax) = O(e^{-C \wb} \Tmax)$. In some cases (problems which are asymptotically complete and have no bound states), we believe that this error can be proved to be time independent.

The main drawback of our algorithm is that it does not provide us the ability to filter waves for which the wavelength is longer than the buffer region. This is a due to the Heisenberg uncertainty principle. This problem will be addressed by a novel multiscale argument in \cite{us:multiscale}. The problem is shared with most methods of absorbing boundaries. The TDPSF algorithm also alerts the user when it fails, unlike all other methods we are aware of.

\section{The Algorithm}

\subsection{The Windowed Fourier Transform}

The discrete windowed Fourier transform frame is a way of localizing a function $f(\vx)$ in the phase space. For the big picture, see \cite{MR1162107,MR924682,MR836025}; see also  \cite[Section 3] {us:TDPSFrigorous} for technical details needed by the TDPSF.

The WFT is an expansion of a function $f(x)$ into the following ``basis'' (actually a frame) for $L^2(\Rn)$:
\begin{eqnarray}
  \label{eq:WFTdef}
  \frmlt{a}{b}{x}=\frmltDef{a}{b}{x},&& \va,\vb \in \ZtwoN
\end{eqnarray}
The lattice spacing in position and momentum, $\xs$, $\ks$, and the standard deviation $\sdev$ must positive real numbers. To be a frame, it must hold that $\xs \ks < 2\pi$. 

This family of functions is not a basis in the usual sense, but a frame. A family of functions $\{ \frmltGen{j}{x} \}_{j \in J}$ is a frame in the Hilbert space $H$ if there exists an $0 < \frmLB \leq \frmUB < \infty$ so that:
\begin{equation*}
  \frmLB \norm{f(x)}{H} \leq \left( \fsumGen{j}{J} \abs{\IP{f}{\frmltGen{j}{x}}}^2 \right)^{1/2} \leq \frmUB \norm{f(x)}{H}
\end{equation*}
A frame is an overcomplete, non-orthonormal basis with a numerically stable reconstruction procedure. We can then decompose any $f \in L^2(\Rn)$ as:
\begin{eqnarray}
  f(\vx) & = & \fsum{a}{b}{\ZtwoN} \fcoeffs{f}{a}{b} \frmlt{a}{b}{x} \\
  \fcoeffs{f}{a}{b} & = &\IP{f(\vx)}{\frmltDual{a}{b}{x}}
\end{eqnarray}
The functions $\frmltDual{a}{b}{x}$ are given by $\frmltDual{a}{b}{x}=\frmltDualDef{a}{b}{x}$
with $\dualW(\vx) \in L^2(\Rn)$ (see \cite{MR1162107,MR924682} for procedures to calculate $\dualW(\vx)$).

In \cite{us:TDPSFrigorous} we extend the analysis of \cite{MR924682} and show that if $\xs \ks = 2\pi/\oversample$ for $\oversample \in 2\mathbb{N}$, $g(\vx)$ decays exponentially in $\absSmallOne{\vx}$ and $\hat{\dualW}(\vk)$ decays exponentially in $\absSmallOne{\vk}$:
\begin{equation}
  \label{eq:dwDecayRate}
  \abs{\partial_x^{\vec{\alpha}} \dualW(\vec{x})} \leq \DWdecayConstXd{\vec{\alpha}} e^{-\DWdecayRateX \absSmallOne{\vx}}
\end{equation}
$\DWdecayConstXd{\vec{\alpha}}$ is a constant (explicitly bounded in \cite{us:TDPSFrigorous}) as is $\DWdecayRateX=\xs \oversample/8\pi \sdev$. \eqref{eq:dwDecayRate} holds for $\hat{\dualW}(\vk)$, but with $\xs$ and $\ks$ swapped and $\sdev^{-1}$ replacing $\sdev$.

\subsubsection{Phase Space Localization}

The WFT allow us to define a concrete realization of phase space. From here onward, we will consider $\ZtwoN$ to be a discrete realization of phase space. The vector $(\vec{a},\vec{b}) \in \ZtwoN$ will represent the point at $\vec{a} \xs$ in position, and $\vec{b} \ks$ in momentum.

With this in mind, we can now construct phase space localization operators very simply. For a set $F \in \ZtwoN$, we define $\psl{F}$ by:
\begin{equation}
  \label{eq:PSLoperatorDef}
  \psl{F} \psi(x) = \fsum{a}{b}{F} \fcoeffs{\psi}{a}{b} \frmlt{a}{b}{x} = \fsum{a}{b}{F} \IP{\frmltDual{a}{b}{x}}{\psi(x)} \frmlt{a}{b}{x} 
\end{equation}

It can be shown \cite{us:TDPSFrigorous,MR966733,MR1162107} that phase space localization in terms of the WFT corresponds closely (though not exactly) to the usual phase space localization. We state one theorem (proved in \cite{us:TDPSFrigorous}, based on a result by Daubechies from \cite{MR966733,MR1162107}) as an example of this:

\newcommand{\Bxp}[0]{B_{X}}
\newcommand{\Bkp}[0]{B_{K}}
\newcommand{\BxpZ}[0]{B_{X'}}
\newcommand{\BkpZ}[0]{B_{K'}}

\begin{theorem}
  \label{thm:pslocWFTbox}
  Let $\Bxp=[-X,X]^N$, $\Bkp=[-K,K]^N$ for $X,K<\infty$. Suppose also that $\xs \ks = 2\pi/\oversample$ for $\oversample \in 2\mathbb{N}$. Then there exists $C$, $\XBuffCube{\epsilon}{K}{X}$ and $\KBuffCube{\epsilon}{K}$ so that if $X'=X-\XBuffCube{\epsilon}{K}{X}$, $K'=K-\KBuffCube{\epsilon}{K}$, then:
    \begin{multline}
    \label{eq:pslocWFTcube}
    \norm{f(x)-\psl{B_{X'} \times B_{K'}} f(x)}{\Hs} \\
    \leq  C \left(\norm{(1-\proj{\Bxp}{\xs}{s}{\vx} )f(\vec{x})}{\Hs} + \norm{(1- \proj{\Bkp}{\ks}{0}{\vk})f(\vec{x})}{\Hs} + \epsilon\norm{f}{\Hs} \right) 
  \end{multline}
  The constant $C$ depends on $\epsilon,s,\xs,\ks,\sdev, X$ and $K$, as do $\XBuffCube{\epsilon}{K}{X}$ and $\KBuffCube{\epsilon}{K}$. $C$ is given explicitly in \cite{us:TDPSFrigorous}.
\end{theorem}

\subsubsection{Computation of the WFT Coefficients and Phase Space Projections: How to do it, and how hard it is}
\label{sec:wftHOWTO}

We now present the algorithm for computing the framelet coefficients. The algorithm consists of computing the products $f(\vx) \dualW(\vx - \va \xs)$, followed by Fourier transforming the results. Due to the spatial decay of $\dualW(\vx)$, we can truncate the domain to the box $[-L_{\epsilon},L_{\epsilon}]^{N}$ with minimal error.

Note that $\dualW(\vx)$ can be calculated as accurately as necessary \cite{MR1162107,MR924682}, and we will not discuss this here.

\begin{algo}
  \label{algo:WFT}
  Calculation of Windowed Fourier Transforms
  
  This algorithm calculates, for a function $f(x)$, the WFT coefficients $\fcoeffs{f}{a}{b}$ for $\va \in A \subseteq \Zn$. We \emph{assume} that the frequencies are bounded above by $\kmax$. The lattice spacing in frequency, $\ks$ is taken to be $2\pi/L_{\epsilon}$; taking it to be any larger than this yields only logarithmic improvements in computational complexity.

  \begin{enumerate}
  \item Let $A \subset \Zn$ be some set of position coordinates.
  \item For each $\va \in A$ multiply $f(\vx)\dualW(\vx-\va \xs)$ for $\vx \in [-L_{\epsilon},L_{\epsilon}]^{N} + \va \xs$ only.
  \item Calculate the Fast Fourier Transform of $f(\vx)\dualW(\vx-\va \xs)$ on this region. For each $\va$, the resulting function is $\fcoeffs{f}{a}{b}$.
  \end{enumerate}
\end{algo}

We observe that this algorithm is local in space. This means that if $A$ is finite, then the computational cost is proportional to:
\begin{equation*}
  \abs{A} \kmax^N \abs{\operatorname{supp} \dualW(\vx) }^N \log (\kmax \abs{\operatorname{supp} \dualW(\vx) }) = \abs{A} 2^{N} \kmax^{N} L_{\epsilon}^{N}
\end{equation*}

The reason for this complexity is as follows. For each $\va \in A$, we need to compute an FFT. The FFT is computed only on the ``support'' (after truncating) of $\dualW(\vx)$ and the lattice spacing in this region is $2\pi/\kmax = O(1/\kmax)$. Thus, there are $O(\abs{\operatorname{supp} \dualW(\vx)}^N \kmax^N)=M$ data points in this region. The FFT has computational complexity $M \log(M)$. In addition, we need to compute $\abs{A}$ of these FFT's.

Note that this algorithm can be parallelized very easily, simply by having different processors compute the FFT's for different $\va \xs$. 

\begin{remark}
\label{rem:complexityOfFilter}
As an example, consider the case when $A=\{\va \in \mathbb{Z}^N : \va \xs \in \FBox \setminus \IBox\}$. Then the size of $A$ is proportional to $\abs{\FBox \setminus \IBox} / \xs^N$. If $\Lb \gg \wb$, then this is of order $\Lb^{N-1}$, and the computational complexity is $O(\Lb^{N-1} \kmax^N \log(\Lb \kmax))$.

We consider this case since this is what the TDPSF requires.
\end{remark}

Phase space projections can also be computed. Let $F \subset \ZtwoN$ be finite. Let $A=\{ \va \in \mathbb{Z}^N : \exists \vb \in \mathbb{Z}, (\va,\vb) \in F \}$. Then we provide the phase space projection algorithm:

\begin{algo}
\label{algo:psl}
  Phase Space Projection Algorithm

  This algorithm computes the phase space projection onto a region of phase space $F$.

  \begin{enumerate}
  \item Compute $A = \{\va : (\va,\vb) \in F \}$. Assume this is finite.
  \item Compute $\fcoeffs{f}{a}{b}$ for $\va \in A$ as in Algorithm \ref{algo:WFT}.
  \item Define a new function $f^{t}: \ZtwoN \rightarrow \mathbb{C}$ by:
    \begin{eqnarray*}
      f^{t}_{(\va,\vb)} = & \fcoeffs{f}{a}{b}, & (\va,\vb) \in F\\
      f^{t}_{(\va,\vb)} = & 0,& (\va,\vb )\nin F
    \end{eqnarray*}
  \item Compute the inverse WFT of $f^{t}$. The result approximates $\psl{F} f(x)$, with errors caused by the truncation of $\dualW(\vx)$ in Algorithm \ref{algo:WFT}.
  \end{enumerate}
\end{algo}

Clearly the computational complexity of Algorithm \ref{algo:psl} is of the same order as that of Algorithm \ref{algo:WFT}.

\subsection{Propagation with Periodic Boundaries: FFT/Split Step Algorithm}
\label{sec:SplitStep}

The TDPSF is a filtering algorithm, which is built on top of another propagation method which solves \eqref{eq:NLSE} on a box with periodic boundaries. The exact manner in which this is done is irrelevant for our purposes, provided it is sufficiently accurate. However, it is important to note that the spectral propagation method is significantly better than the usual FDTD methods, and that compatability with spectral methods is an important feature of the TDPSF.

Fix a grid spacing $\delta x$ and timestep $\delta t$. The computational grid is simply the region $\CBox$ sampled at the lattice with spacing $\delta x$. This corresponds to a lattice spacing in momentum of $2\pi/\LbC$, with maximal momentum $2\pi / \delta x$. A common rule of thumb is that if the problem has a maximal momentum $\kmax$, then $\delta x = 4 \pi/ \kmax$.

\begin{algo} Split Step Propagation Algorithm
  \label{algo:splitStep}
  
  This algorithm approximates the propagation of \eqref{eq:NLSE} on the region\\ $\CBox$. The accuracy is $O(\delta t^{2})$ time and $O(\delta x^{\omega})$ (``spectral accuracy'') in space. It is $O(\delta t^{3})$ in time if \eqref{eq:NLSE} is linear.

  First, the discrete operator $e^{i \Lap t} f(x)$ is defined by computing the FFT of $f(x)$, multiplying by $e^{i k^{2} t}$ and then computing the inverse FFT of the result.

  With this in mind, the algorithm is as follows (taking $\psi_{0}(x)$ as initial condition, and assuming $t$ is a multiple of $\delta t$):
  \begin{enumerate}
  \item Define $\psi_{1/2}(x)=\Uf{\delta t/2} \psi_{0}(x)$.
  \item For $j=0,\ldots,t/\delta t-2$, define:
    \begin{equation*}
      \psi_{j+1+1/2}(x) = \Uf{\delta t} e^{-i \nlinsemi{j \delta t}{\psi_{j+1/2}(x)} \delta t} \psi_{j+1/2}(x)
    \end{equation*}
  \item Finally, define $\psi_{t/\delta t}(x) \approx \psi(x,t)$ by:
    \begin{equation*}
      \psi_{t/\delta t}(x) = \Uf{(\delta t/2)} e^{-i \nlinsemi{j \delta t}{\psi_{j+1/2}(x)} \delta t} \Uf{\delta t} \psi_{t/\delta t-3/2}(x)
    \end{equation*}
  \end{enumerate}
  This is a numerical realization of the Trotter-Kato product formula:
  \begin{equation}
    \label{eq:trotterKatoProductFormula}
    \U{t} \approx \Uf{(-\delta t/2)}\left[
        \prod_{j=1}^{t/\delta t} \Uf{\delta t} e^{-i \nlinsemi{j \delta t}{\psi(x,j \delta t)} \delta t}
      \right] \Uf{\delta t/2}
    \end{equation}
\end{algo}

The computational complexity of Algorithm \ref{algo:splitStep} is
\begin{equation*}
  O((\Lb+\wb)^{N} \kmax^{N} \ln[(\Lb+\wb)\kmax] ) 
\end{equation*}
per timestep (since the number of data points is of order $(\Lb+\wb)\kmax$). The number of timesteps is $\Tmax/\delta t$. Thus, split step propagation has time complexity $O[\Lb^N \kmax^N \log(\Lb \kmax) (\Tmax/\delta t)]$. Algorithm \ref{algo:splitStep} is by now a textbook result \cite{boyd:spectralmethods}, so we will not discuss it further. 

\subsection{The TDPSF Algorithm}

The TDPSF algorithm consists of solving \eqref{eq:NLSE} on the box $\CBox$ with periodic boundary conditions, using Algorithm \ref{algo:splitStep} or any other method. At times $\Tstep$, $2 \Tstep$, etc, we subtract from the solution the projection onto Gaussians which are outgoing by means of Algorithm \ref{algo:psl}.

Outgoing Gaussians are defined to be those $\frmlt{a}{b}{x}$ with $\va \xs \in \FBox \setminus \IBox$ and which are moving strictly outward under the free flow; that is, those which are leaving $\IBox$ soon (before a time $\Tstep$ has passed) and which will never return.

Incoming Gaussians are those which will not leave $\CBox$ before time $\Tstep$ has passed.

There is also a set of ambiguous gaussians. These are gaussians with velocity sufficiently low so that they spread about their center faster than they move. They are ambiguous because they are both leaving $\CBox$ and returning to $\IBox$ due to their spread. An example of this would be $e^{-(x-\Lb-\wb/2)^{2}/2\sdev^{2}}$. This gaussian is located inside the buffer region, and is spreading laterally with no motion whatsoever. 

The TDPSF algorithm does not know what to do with these gaussians, which correspond to low frequency waves. However, unlike other algorithms (e.g. Dirichlet-Neumann, PML or absorbing potential), the TDPSF does know when such waves are causing a problem. Thus, incorrect results are never returned.

\begin{algo} Propagation algorithm
  \label{algo:TDPSFpropagator}
  This is the main propagation algorithm. In this algorithm, we consider $\sdev$, $\xs$, $\ks$ to be fixed. $\Lb$ and $\Tmax$ are also considered fixed. The initial data is considered fixed, and localized inside $\IBox$. The approximation will be denoted by $\Psi(x,t)$. Also, fix a small tolerance $\epsilon > 0$.
  
  \begin{enumerate}
  \item Before beginning, precalculate the set of framelets which are outgoing, and those which are ambiguous.
  \item Define $\Psi(x,t)$ iteratively as follows. Loop over $n = 0\ldots \Tmax/\Tstep$. In what follows, the propagator $\U{t}$ is calculated by Algorithm \ref{algo:splitStep}.
    \begin{enumerate}
    \item For $t \in [(n-1)\Tstep,n\Tstep)$, define 
      \begin{equation*}
        \psi(x,t) = \U{t-(n-1)\Tstep} \psi(x,(n-1)\Tstep)
      \end{equation*}

    \item For $t = n \Tstep$, define:
      \begin{equation*}
        \Psi(x,n \Tstep) = (1 - \psl{\operatorname{OUT}}) \U{\Tstep}\psi(x,(n-1)\Tstep)
      \end{equation*}

      \item At each integer multiple of $\Tstep$, compute 
        \begin{equation*}
          \norm{\psl{\operatorname{AMB}} \U{\Tstep}\psi(x,(n-1)\Tstep)}{\Hs}
        \end{equation*}
        where $\operatorname{AMB}$ is the set of ambiguous framelets. If this quantity is greater than $\epsilon$, \emph{stop the program and notify the user}.
    \end{enumerate}
  \end{enumerate}
\end{algo}

\begin{remark}
  The TDPSF fails gracefully in the following respect. If the TDPSF is unable to filter due to slowly moving waves, the program fails and the user is notified. Absorbing potentials, the PML and Dirichlet-Neumann maps do not have this property.
\end{remark}

\begin{remark}
  In order to calculate $\psl{\operatorname{OUT}}$, we must calculate the WFT coefficients of the solution in the region $\CBox \setminus \IBox$. To calculate $\psl{\operatorname{AMB}}$, we need the WFT coefficients from the same region. Thus, only one WFT is needed for steps (b) and (c). The cost of computing the $H^{s}$ norm of the ambiguous framelets is cheaper than the cost of a WFT; this is merely a sum rather than a set of FFT's.
\end{remark}

Thus, all that remains is to do is explain which framelets are outgoing, incoming and ambiguous.

For a given $(\va,\vb)$, this can be determined easily, by a direct calculation based on \eqref{eq:freeFrameletProp}. This is done in \cite{us:TDPSFrigorous} for $s=0,1$, and other cases can be done with the help of Maple/Mathematica.

The end result of the calculation is the following. Define $R_{\vb}(t)$ by
\begin{subequations}
  \begin{equation}
    R_{\vb}(t) = \sqrt{\sdev^2+t^2/\sdev^2}(\Gammai(N/2,2 \epsilon^2 \pi^{N/2} / \abs{S^{N-1}}) )^{1/2}
  \end{equation}
  if we measure the error in $L^{2}(\Rn)$ or
  \begin{multline}
    R_{\vb}(t) = \sqrt{\sdev^2+t^2/\sdev^2} \max\Bigg\{
    \left[\Gammai\left(N/2,  \frac{\epsilon^2 \pi^N/2}{2 \abs{S^{N-1}}(1+\absSmallTwo{\vb \ks}^2)}\right)\right]^{1/2}, \\
    \left[\Gammai\left((N+2)/2, \frac{\epsilon^2 \sdev^2 \pi^{N/2}}{2 \abs{S^{N-1}}}\right)\right]^{1/2}
    \Bigg\} 
  \end{multline}
\end{subequations}
if we measure the error in $H^{1}(\Rn)$.

The function $\Gammai(a,x)$ is the inverse function of $x \mapsto \Gamma(a,x)$, with $\Gamma(a,x)$ the complementary incomplete Gamma function. This definition implies that all the mass (up to a tolerance $\epsilon$) of $\Uf{t} \frmlt{a}{b}{x}$ is, at time $t$, contained in a ball of radius $R_{\vb}$ about the point $\va \xs + \vb \ks t$. That is:
\begin{subequations}
  \label{eq:freeFrameletPropBB}
  \begin{equation}
    \norm{\Uf{t} \frmltDef{a}{b}{x} }{H^{0,1}({B(t)^C})} \leq \epsilon
  \end{equation}
  \begin{equation}
    B(t) = \{x : \absSmallTwo{\va \xs + \vb \ks t - \vx} \leq R_{\vb}(t) \}
  \end{equation}
\end{subequations}
Thus, if $B(t)$ does not intersect $\IBox$, this framelet is strictly bad. 

Similar calculations can be done in $\Hs$, although $R_{\vb}$ must be redefined (and now depends on $\vb$). Thus, given the result of this calculation, we define the set of outgoing Gaussian's  to be the set of $(\va,\vb) \in \ZtwoN$ such that:

\begin{enumerate}
\item $\va \xs \in \FBox \setminus \IBox$
\item $B(t)$ does not intersect $\IBox$ for any $t < \Tmax$.
\end{enumerate}

A simple criterion to determine whether a given framelet is outgoing is to check whether $\va \xs \in \FBox \setminus \IBox$, and whether $\vb_{n} \ks \geq 2 \ln \epsilon/\sigma$. Here, $\vb_{n}$ is the component of $\vb \ks$ pointing in the direction normal to $\IBox$. Conversely, this means that to filter waves with frequency $\kmin$ (or higher) making an error $\epsilon$, then $\sdev \geq \ln \epsilon / \kmin$. This determines the width of the buffer region. This appears to be somewhat larger than what a PML requires, although we were unable to obtain high accuracy using current implementations of the PML (see Remark \ref{remark:commentsOnPML}).

The set of ambiguous framelets is the set of framelets for which $B(t)$ intersects $\IBox$, but also intersects $\Rn \setminus \FBox$. That is to say, these are framelets with both incoming and outgoing components.

A calculation based on \eqref{eq:freeFrameletProp} shows that these consist of framelets with either low velocities, or framelets which are located on the boundary but are moving tangentially to the boundary.

To assess the time complexity of the TDPSF algorithm, we recall remark \ref{rem:complexityOfFilter}, which says that when $\Lb \gg \wb$, the cost of computing $\psl{\operatorname{OUT}}$ (by Algorithm \ref{algo:psl}) is only $O(\Lb^{N-1} \kmax^{N} \log(\Lb \kmax))$. Thus, in this case, we find that the added complexity of TDPSF propagation over the Split Step algorithm is of lower order than the Split Step itself. However the constant is significantly larger in our experiments, so on a small grid the TDPSF propagator is slower than FFT/Split Step propagation.

We discuss some possible improvements on this algorithm in section \ref{sec:outlook}.

\section{Numerical Examples}
\label{sec:tests}
In this section we discuss the results of our numerical tests.

The TDPSF algorithm is built in the program Kitty. Kitty is written in Python, with C extensions, calling the external libraries FFTW \cite{FFTW05}, Numarray and Matplotlib. The external programs gnuplot, ImageMagick and gifsicle were also used for making graphs/movies.

Kitty is licensed under the GPL. It is very much a work in progress, and has little documentation and minimal user interface (an end-user version is currently in progress). Various test cases, spanning many types of parameters, are also available for download from the author's webpage, http://math.rutgers.edu/\~{}stucchio. 

\subsection{Simple Tests: Free Schr\"odinger Equation}
\label{sec:testsSimpleTR}

The standard method for testing absorbing boundaries is simply to throw coherent states (which are well localized in frequency) at the boundary and compute the difference between the approximate result and the true solution. This is a useful test, although it by no means completely characterizes the errors, as we discuss in Section \ref{sec:testsMediumRange}.

Our implementation is as follows. The free Schr\"odinger equation (i.e.\\ $\nlinsemi{t}{\psi(x,t)}=0$) was solved on $[-102.4,102.4]$  with TDPSF boundaries on the regions $[-102.4,-88]$ and $[88,102.4]$. The lattice spacing was $\delta x=0.1$ in space, $5 \times 10^{-4}$ in time, and the TDPSF filtering interval was $\Tstep=2 \times 10^{-3}$. The initial condition was taken to be $\psi(x,0)=e^{-x^2/4} e^{i v x}$ with $v$ ranging from $1$ to $25$. The simulation was run from $t=0$ out to $\Tmax = 3 \cdot 51.2/v$.

\begin{figure}
\setlength{\unitlength}{0.240900pt}
\ifx\plotpoint\undefined\newsavebox{\plotpoint}\fi
\sbox{\plotpoint}{\rule[-0.200pt]{0.400pt}{0.400pt}}%
\input{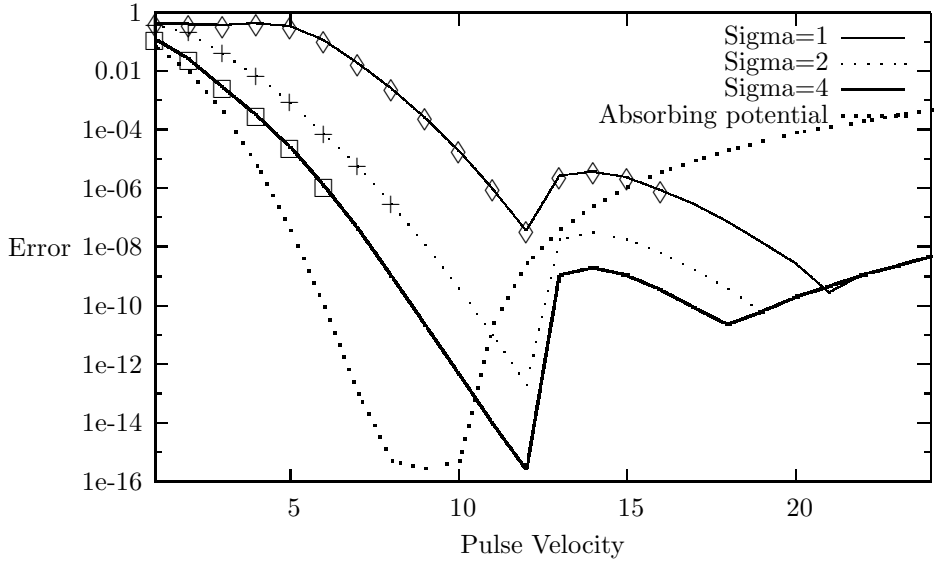}
\caption{A graph of error vs the velocity of an outgoing pulse. The boxes, diamonds and ``+'' signs indicate TDPSF error reports.
}
\label{fig:tPlusRTest}
\end{figure}

After the solution was given sufficient time to exit the computational domain, the error in the region $[-88,88]$ was measured (comparing with the exact result \eqref{eq:freeFrameletPropBB}). The error as a function of velocity is graphed in Figure \ref{fig:tPlusRTest}. The TDPSF was used with $\sigma=1,2,4$. The results are comparable to the complex absorbing potential $V(x)=-25i e^{-(x-25.6)^2/16}$, which is also shown in Figure \ref{fig:tPlusRTest}. The width of the complex potential was chosen so that it's spatial extent is comparable to the width of the TDPSF used. Altering the height of the potential changed only position of the dip in Figure \ref{fig:tPlusRTest}, and achieved little other benefit.

In step 2c of Algorithm \ref{algo:TDPSFpropagator}, the TDPSF checks whether it is making an error. The above simulations were run with tolerance $10^{-6}$, and the simulations for which the TDPSF reported an error are indicated in Figure \ref{fig:tPlusRTest}. Each time the TDPSF failed (and 4 times when it was successful), it notified the user of the error. The absorbing potential was unaware of it's own failures.

This particular example demonstrates no major advantage of the TDPSF over the absorbing potential apart from the awareness of errors when the occur. The TDPSF works better for high frequencies, but not low. The advantage of the TDPSF is not that it succesfully dissipates outgoing waves, but that it does not dissipate incoming waves. An example of this will be demonstrated in the next section. 

\begin{remark}
  \label{remark:commentsOnPML}
  Using the programs from \cite{MR2032866} (available at http://www.math.unm.edu/$\sim$hagstrom/), we ran a similar set of simulations to the ones above using the PML for boundary absorption. Since FDTD methods rather than spectral methods were used, we took $\delta x = 6.7 \times 10^{-4}$, $\delta t=2.4 \times 10^{-6}$, and a PML with a thickness of 500 lattice points\footnote{The computational region was $[-20,20]$ with $30,000$ lattice points. Solving on $[-100,100]$ would require require $5 \times 30,000=150,000$ lattice points (the TDPSF/spectral solver uses $512$ lattice points for the region $[-20,20]$). }. The error was of order $10^{-2}-10^{-3}$, probably due to the use of FDTD instead of spectral methods for solving the interior problem. 
\end{remark}

\subsection{Harder Tests: Long Range Potentials}
\label{sec:testsMediumRange}

The main problem associated with an absorbing potential or PML is that not all waves located near the boundary are outgoing. The problem is that some waves are incoming, and should not be dissipated. 

Consider the following linear Schr\"odinger equation (with $(\vx,t) \in \mathbb{R}^{2+1}$):
\begin{eqnarray}
  i \partial_t \psi(x,t)& =& \left[-\Lap - \frac{15}{0.05 \absSmallTwo{\vx}^2+1} \right] \psi(x,t) \label{eq:testMediumRange} \\
  \psi(x,0)& =& e^{i 7 x_2 } e^{-\absSmallTwo{\vx}^2/20} + e^{i 4 x_1} e^{-\absSmallTwo{\vx}^2/20}\nonumber
\end{eqnarray}

The initial condition consists of two coherent states of equal mass, one with velocity $4$ and one with velocity $7$. The notable fact about this particular potential is that the fast Gaussian has enough kinetic energy to (mostly) escape from the binding potential, while the slow Gaussian does not. The slow Gaussian moves toward the boundary, turns around and returns.

The problem with the absorbing potential approach is that it does not distinguish between incoming and outgoing waves. It dissipates everything on the boundary including the waves that should have returned. We believe this should also occur for the PML (at least the PML found in \cite{MR2032866}).

We ran three simulations of \eqref{eq:testMediumRange}. The first was performed using the TDPSF with $\sdev=2.0$, $\xs=0.8$, $\ks = 2\pi/12.8=0.491$ and $\Tstep=0.1$. The region of computation was $[-25.6,25.6]^2$, with lattice spacing $\delta x = 0.2$ in space and $\delta t = 0.025$ in time. The region of interest was $[-10,10]^{2}$. Note that even if the region of interest were larger (as in Section \ref{sec:testsSimpleTR}), the width of the TDPSF region would not grow\footnote{We could have performed this test on the box $[-204.8,204.8]^{2}$, with the TDPSF region still having width $15.6$. But comparisons with an accurate simulation on a larger box would be impossible.}. An identical simulation was performed (on the same region) with an absorbing potential 
\begin{equation*}
  V_1(\vx) = -20i e^{-(\vx_1 \pm 25.6)^2/36} -20i e^{-(\vx_2 \pm 25.6)^2/36}
\end{equation*}

The third was solved with periodic boundary conditions on $[-204.8,204.8]^2$ (all other parameters were the same). This boundary is sufficiently distant so that the outgoing waves cannot return to the origin for a time $408.6/7.0 \approx 58$. Letting $\psi_{t}(x,t)$ be the solution with the TDPSF boundary, $\psi_{a}(x,t)$ be the solution with the absorbing boundary and $\psi_{d}(x,t)$ be the solution with the distant boundary, we measured the relative error:
\begin{equation*}
  E_{t,a}(t) = \frac{
    \norm{\psi_{t,a}(x,t) - \psi_{d}(x,t)}{L^{2}([-10,10]^{2})}
    }{
      \norm{\psi_{0}(x)}{L^{2}(\mathbb{R}^{2})}
  }
\end{equation*}

\begin{figure}
\setlength{\unitlength}{0.240900pt}
\ifx\plotpoint\undefined\newsavebox{\plotpoint}\fi
\sbox{\plotpoint}{\rule[-0.200pt]{0.400pt}{0.400pt}}%
\input{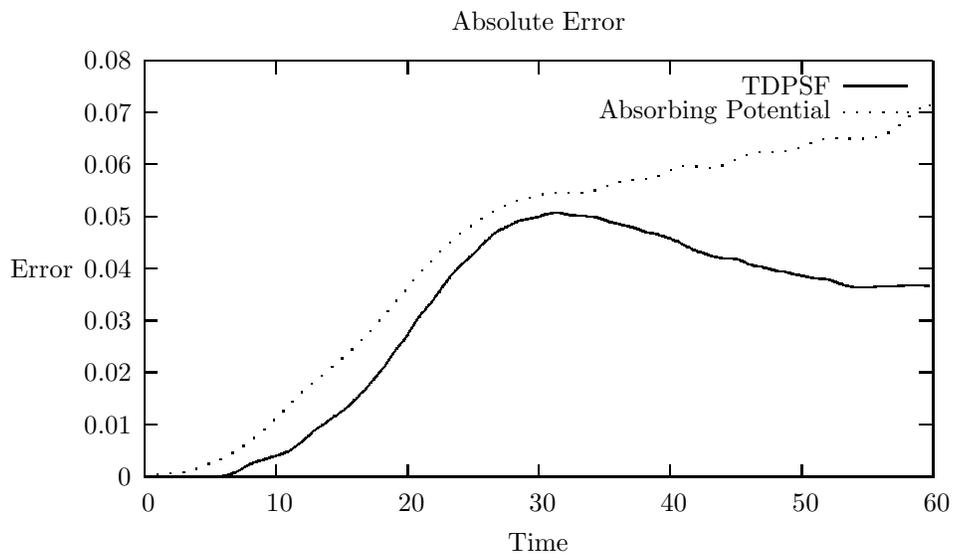}
\caption{
  A graph of the error in $L^{2}([-10,10]^{2})$ for the absorbing potential and TDPSF up to $t=58$.
}
\label{fig:mediumRangeError}
\end{figure}

Figure \ref{fig:mediumRangeError} plots the error. In total, the error for the absorbing potential reaches $7.0 \%$ by $t=58$, while the TDPSF makes an error of $5.0 \%$. The reason the error is so large for the TDPSF is that there is a concentration of mass near $k=0$ near the boundary. This poses a problem for the absorbing potential, as Figure \ref{fig:mediumRangeError} shows.

However, the problem with the absorbing potential is even more serious than it appears at first glance. The true solution $\psi(x,t)$ can be written as $\psi(x,t)=\sum_{j} e^{-i E_{j} t}\psi_{j} \phi_{j}(x) + \psi_{d}(x,t)$, with $\phi_{j}(x)$ the bound states and $\psi_{d}(x,t)$ the dispersive part of $\psi(x,t)$. According to our simulations, by approximately $t=30$, $\psi_{d}(x,t)$ appears to have dispersed. Therefore, after this point, all that should remain are bound states.

\begin{figure}
\setlength{\unitlength}{0.240900pt}
\ifx\plotpoint\undefined\newsavebox{\plotpoint}\fi
\sbox{\plotpoint}{\rule[-0.200pt]{0.400pt}{0.400pt}}%
\input{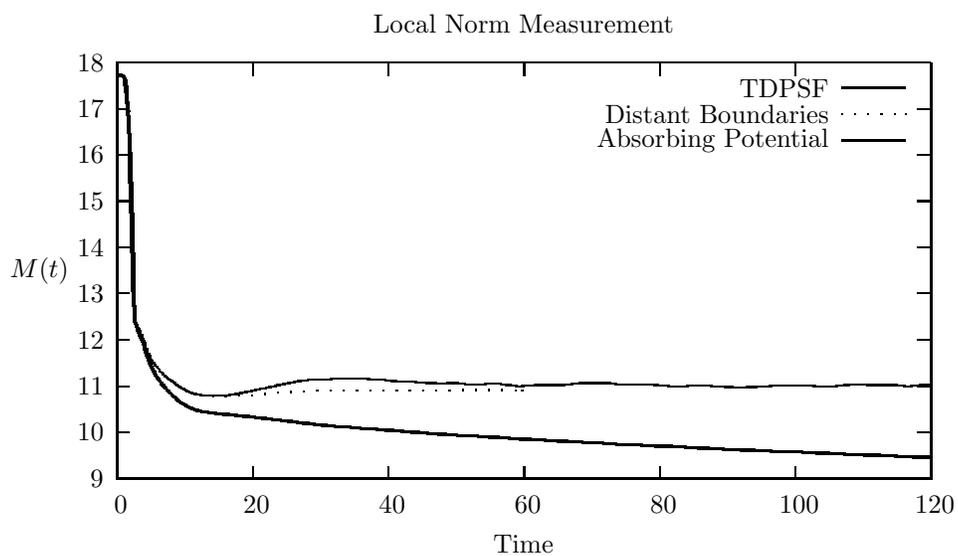}
\caption{
  A graph of $M(t)=\norm{\psi(x,t)}{L^2([-10,10]^2)}$ vs $t$. The distant boundary simulation is invalid at time $t=58$, due to the return of outgoing pulse.
}
\label{fig:mediumRangeL2loc}
\end{figure}

This implies that $M(t)=\norm{\psi(x,t)}{L^2([-10,10]^2)}$ should remain approximately constant after $t=30$. This behavior can be seen in the distant boundaries simulation as well as the TDPSF simulation. This does not occur for the absorbing potential. Figure \ref{fig:mediumRangeL2loc} graphs $M(t)$ for $t \in [0,60]$ (in the distant boundaries case) and $t \in [0,120]$ otherwise. The absorbing potential appears to be dissipating, although we know this is impossible. For larger times, even the qualitative behavior of the absorbing potential simulation is wrong.

The reason the TDPSF performs better than the complex potential is that it distinguishes outgoing waves from incoming waves. The TDPSF only removes waves which sit on the boundary and are also clearly outgoing. Low frequency waves on the boundary, which return to the computation region are incorrectly dissipated by the absorbing potential, but (correctly) ignored by the TDPSF. This is why the TDPSF gives the correct long time behavior.

This problem is more dangerous than this example suggests. Although in this example we can determine that the dissipation is artificial, we cannot always do so. Similar problems occur often in atomic physics, but with exponentially decaying resonances as well as bound states (corresponding to complex $E_{j}$). In this case, one cannot distinguish the dissipation caused by the absorbing potential from that caused by the dynamics, and the measured decay rate will be incorrect. Since the object of such a simulation is often the measurement of the decay rate, this is a serious problem.

\subsection{Soliton Filtering}
\label{sec:testsSoliton}

Consider the nonlinearity, $\nlinsemi{t}{\psi(\vx,t)}=-\abs{\psi(\vx,t)}^2$. It is desirable to construct a numerical algorithm which filters outgoing solitons as well as free waves. Although the TDPSF was not designed to filter solitons, it turns out to work well for solitons which are not moving slowly. This is however a coincidence, and cannot be expected to hold in general. A referee pointed out that it will fail for the KdV equation.

The reason for this is that an outgoing soliton with sufficiently high velocity is localized in phase space on outgoing waves. Consider a simple soliton, $\phi(x,v,t) = 2^{-1/2} e^{i(v x+(1-v^2)t)}/\cosh((x-vt))$. A simple calculation shows that for $\abs{k-v} \gg 1$, $\hat{\phi}(k,v,t) \sim e^{-\abs{k-v}}$. Thus, for $k \gg \kmin$, this shows that the framelet coefficients of $\phi(x,v,t)$ which are moving too slowly to resolve have exponentially small mass. This shows that under the free flow this soliton is strictly outgoing.

The soliton is also leaving the box under the full flow $\U{t}$. Although $\Uf{t}$ and $\U{t}$ move the soliton very differently (one dispersively, one coherently), they both move it out of the box and in nearly the same direction. For this reason, we expect the TDPSF is to filter soliton solutions correctly.

We ran numerical tests to demonstrate this as follows. We solved \eqref{eq:NLSE} with $\nlinsemi{t}{\psi(\vx,t)}=-\abs{\psi(\vx,t)}^2$ on the region $[-25.6,25.6]$ with $\delta x = 0.05$, $\delta t=0.002/v$ and $\Tstep=0.008/v$ (the timestep's are scaled with the velocity to speed up the simulations). In this simulation, $\Lb=12.0$ and $\wb=13.6$. The initial condition was taken to be $\psi(x,0)=2^{-1/2} e^{i v x} / \cosh(x)$ for $v=1..15$. 

\begin{figure}
\label{fig:solitonAbsorption}
\setlength{\unitlength}{0.240900pt}
\ifx\plotpoint\undefined\newsavebox{\plotpoint}\fi
\sbox{\plotpoint}{\rule[-0.200pt]{0.400pt}{0.400pt}}%
\setlength{\unitlength}{0.240900pt}
\ifx\plotpoint\undefined\newsavebox{\plotpoint}\fi
\begin{picture}(1500,900)(0,0)
\sbox{\plotpoint}{\rule[-0.200pt]{0.400pt}{0.400pt}}%
\put(221.0,123.0){\rule[-0.200pt]{4.818pt}{0.400pt}}
\put(201,123){\makebox(0,0)[r]{ 1e-10}}
\put(1419.0,123.0){\rule[-0.200pt]{4.818pt}{0.400pt}}
\put(221.0,188.0){\rule[-0.200pt]{4.818pt}{0.400pt}}
\put(201,188){\makebox(0,0)[r]{ 1e-09}}
\put(1419.0,188.0){\rule[-0.200pt]{4.818pt}{0.400pt}}
\put(221.0,254.0){\rule[-0.200pt]{4.818pt}{0.400pt}}
\put(201,254){\makebox(0,0)[r]{ 1e-08}}
\put(1419.0,254.0){\rule[-0.200pt]{4.818pt}{0.400pt}}
\put(221.0,319.0){\rule[-0.200pt]{4.818pt}{0.400pt}}
\put(201,319){\makebox(0,0)[r]{ 1e-07}}
\put(1419.0,319.0){\rule[-0.200pt]{4.818pt}{0.400pt}}
\put(221.0,385.0){\rule[-0.200pt]{4.818pt}{0.400pt}}
\put(201,385){\makebox(0,0)[r]{ 1e-06}}
\put(1419.0,385.0){\rule[-0.200pt]{4.818pt}{0.400pt}}
\put(221.0,450.0){\rule[-0.200pt]{4.818pt}{0.400pt}}
\put(201,450){\makebox(0,0)[r]{ 1e-05}}
\put(1419.0,450.0){\rule[-0.200pt]{4.818pt}{0.400pt}}
\put(221.0,515.0){\rule[-0.200pt]{4.818pt}{0.400pt}}
\put(201,515){\makebox(0,0)[r]{ 1e-04}}
\put(1419.0,515.0){\rule[-0.200pt]{4.818pt}{0.400pt}}
\put(221.0,581.0){\rule[-0.200pt]{4.818pt}{0.400pt}}
\put(201,581){\makebox(0,0)[r]{ 0.001}}
\put(1419.0,581.0){\rule[-0.200pt]{4.818pt}{0.400pt}}
\put(221.0,646.0){\rule[-0.200pt]{4.818pt}{0.400pt}}
\put(201,646){\makebox(0,0)[r]{ 0.01}}
\put(1419.0,646.0){\rule[-0.200pt]{4.818pt}{0.400pt}}
\put(221.0,712.0){\rule[-0.200pt]{4.818pt}{0.400pt}}
\put(201,712){\makebox(0,0)[r]{ 0.1}}
\put(1419.0,712.0){\rule[-0.200pt]{4.818pt}{0.400pt}}
\put(221.0,777.0){\rule[-0.200pt]{4.818pt}{0.400pt}}
\put(201,777){\makebox(0,0)[r]{ 1}}
\put(1419.0,777.0){\rule[-0.200pt]{4.818pt}{0.400pt}}
\put(308.0,123.0){\rule[-0.200pt]{0.400pt}{4.818pt}}
\put(308,82){\makebox(0,0){ 2}}
\put(308.0,757.0){\rule[-0.200pt]{0.400pt}{4.818pt}}
\put(482.0,123.0){\rule[-0.200pt]{0.400pt}{4.818pt}}
\put(482,82){\makebox(0,0){ 4}}
\put(482.0,757.0){\rule[-0.200pt]{0.400pt}{4.818pt}}
\put(656.0,123.0){\rule[-0.200pt]{0.400pt}{4.818pt}}
\put(656,82){\makebox(0,0){ 6}}
\put(656.0,757.0){\rule[-0.200pt]{0.400pt}{4.818pt}}
\put(830.0,123.0){\rule[-0.200pt]{0.400pt}{4.818pt}}
\put(830,82){\makebox(0,0){ 8}}
\put(830.0,757.0){\rule[-0.200pt]{0.400pt}{4.818pt}}
\put(1004.0,123.0){\rule[-0.200pt]{0.400pt}{4.818pt}}
\put(1004,82){\makebox(0,0){ 10}}
\put(1004.0,757.0){\rule[-0.200pt]{0.400pt}{4.818pt}}
\put(1178.0,123.0){\rule[-0.200pt]{0.400pt}{4.818pt}}
\put(1178,82){\makebox(0,0){ 12}}
\put(1178.0,757.0){\rule[-0.200pt]{0.400pt}{4.818pt}}
\put(1352.0,123.0){\rule[-0.200pt]{0.400pt}{4.818pt}}
\put(1352,82){\makebox(0,0){ 14}}
\put(1352.0,757.0){\rule[-0.200pt]{0.400pt}{4.818pt}}
\put(221.0,123.0){\rule[-0.200pt]{293.416pt}{0.400pt}}
\put(1439.0,123.0){\rule[-0.200pt]{0.400pt}{157.549pt}}
\put(221.0,777.0){\rule[-0.200pt]{293.416pt}{0.400pt}}
\put(221.0,123.0){\rule[-0.200pt]{0.400pt}{157.549pt}}
\put(40,450){\makebox(0,0){E(v)}}
\put(830,21){\makebox(0,0){Velocity}}
\put(830,839){\makebox(0,0){Mass not absorbed}}
\put(1279,737){\makebox(0,0)[r]{sigma=1}}
\put(1299.0,737.0){\rule[-0.200pt]{24.090pt}{0.400pt}}
\put(1439,138){\usebox{\plotpoint}}
\multiput(1434.68,138.58)(-1.181,0.498){71}{\rule{1.041pt}{0.120pt}}
\multiput(1436.84,137.17)(-84.840,37.000){2}{\rule{0.520pt}{0.400pt}}
\multiput(1349.63,175.58)(-0.588,0.499){145}{\rule{0.570pt}{0.120pt}}
\multiput(1350.82,174.17)(-85.816,74.000){2}{\rule{0.285pt}{0.400pt}}
\multiput(1263.92,249.00)(-0.499,0.517){171}{\rule{0.120pt}{0.514pt}}
\multiput(1264.17,249.00)(-87.000,88.934){2}{\rule{0.400pt}{0.257pt}}
\multiput(1175.91,339.58)(-0.505,0.499){169}{\rule{0.505pt}{0.120pt}}
\multiput(1176.95,338.17)(-85.953,86.000){2}{\rule{0.252pt}{0.400pt}}
\multiput(1088.80,425.58)(-0.537,0.499){159}{\rule{0.530pt}{0.120pt}}
\multiput(1089.90,424.17)(-85.901,81.000){2}{\rule{0.265pt}{0.400pt}}
\multiput(1001.61,506.58)(-0.596,0.499){143}{\rule{0.577pt}{0.120pt}}
\multiput(1002.80,505.17)(-85.803,73.000){2}{\rule{0.288pt}{0.400pt}}
\multiput(914.33,579.58)(-0.680,0.499){125}{\rule{0.644pt}{0.120pt}}
\multiput(915.66,578.17)(-85.664,64.000){2}{\rule{0.322pt}{0.400pt}}
\multiput(826.86,643.58)(-0.822,0.498){103}{\rule{0.757pt}{0.120pt}}
\multiput(828.43,642.17)(-85.430,53.000){2}{\rule{0.378pt}{0.400pt}}
\multiput(739.30,696.58)(-0.992,0.498){85}{\rule{0.891pt}{0.120pt}}
\multiput(741.15,695.17)(-85.151,44.000){2}{\rule{0.445pt}{0.400pt}}
\multiput(650.23,740.58)(-1.625,0.497){51}{\rule{1.389pt}{0.120pt}}
\multiput(653.12,739.17)(-84.117,27.000){2}{\rule{0.694pt}{0.400pt}}
\multiput(539.69,767.59)(-9.616,0.477){7}{\rule{7.060pt}{0.115pt}}
\multiput(554.35,766.17)(-72.347,5.000){2}{\rule{3.530pt}{0.400pt}}
\put(395,770.17){\rule{17.500pt}{0.400pt}}
\multiput(445.68,771.17)(-50.678,-2.000){2}{\rule{8.750pt}{0.400pt}}
\put(308,768.67){\rule{20.958pt}{0.400pt}}
\multiput(351.50,769.17)(-43.500,-1.000){2}{\rule{10.479pt}{0.400pt}}
\put(221,769.17){\rule{17.500pt}{0.400pt}}
\multiput(271.68,768.17)(-50.678,2.000){2}{\rule{8.750pt}{0.400pt}}
\put(1279,696){\makebox(0,0)[r]{sigma=2}}
\multiput(1299,696)(20.756,0.000){5}{\usebox{\plotpoint}}
\put(1399,696){\usebox{\plotpoint}}
\put(1439,133){\usebox{\plotpoint}}
\multiput(1439,133)(-20.734,0.953){5}{\usebox{\plotpoint}}
\multiput(1352,137)(-20.734,0.953){4}{\usebox{\plotpoint}}
\multiput(1265,141)(-20.721,1.191){4}{\usebox{\plotpoint}}
\multiput(1178,146)(-20.228,4.650){4}{\usebox{\plotpoint}}
\multiput(1091,166)(-17.361,11.375){5}{\usebox{\plotpoint}}
\multiput(1004,223)(-15.454,13.855){6}{\usebox{\plotpoint}}
\multiput(917,301)(-14.846,14.505){6}{\usebox{\plotpoint}}
\multiput(830,386)(-15.104,14.236){6}{\usebox{\plotpoint}}
\multiput(743,468)(-15.810,13.448){5}{\usebox{\plotpoint}}
\multiput(656,542)(-16.536,12.544){5}{\usebox{\plotpoint}}
\multiput(569,608)(-17.453,11.234){5}{\usebox{\plotpoint}}
\multiput(482,664)(-18.349,9.702){5}{\usebox{\plotpoint}}
\multiput(395,710)(-18.940,8.490){5}{\usebox{\plotpoint}}
\multiput(308,749)(-20.122,5.088){4}{\usebox{\plotpoint}}
\put(221,771){\usebox{\plotpoint}}
\sbox{\plotpoint}{\rule[-0.400pt]{0.800pt}{0.800pt}}%
\sbox{\plotpoint}{\rule[-0.200pt]{0.400pt}{0.400pt}}%
\put(1279,655){\makebox(0,0)[r]{sigma=3}}
\sbox{\plotpoint}{\rule[-0.400pt]{0.800pt}{0.800pt}}%
\put(1299.0,655.0){\rule[-0.400pt]{24.090pt}{0.800pt}}
\put(1439,133){\usebox{\plotpoint}}
\put(1352,133.34){\rule{17.600pt}{0.800pt}}
\multiput(1402.47,131.34)(-50.470,4.000){2}{\rule{8.800pt}{0.800pt}}
\multiput(1293.39,138.38)(-14.193,0.560){3}{\rule{14.120pt}{0.135pt}}
\multiput(1322.69,135.34)(-57.693,5.000){2}{\rule{7.060pt}{0.800pt}}
\multiput(1206.39,143.38)(-14.193,0.560){3}{\rule{14.120pt}{0.135pt}}
\multiput(1235.69,140.34)(-57.693,5.000){2}{\rule{7.060pt}{0.800pt}}
\multiput(1129.02,148.39)(-9.504,0.536){5}{\rule{11.800pt}{0.129pt}}
\multiput(1153.51,145.34)(-62.509,6.000){2}{\rule{5.900pt}{0.800pt}}
\multiput(1042.02,154.39)(-9.504,0.536){5}{\rule{11.800pt}{0.129pt}}
\multiput(1066.51,151.34)(-62.509,6.000){2}{\rule{5.900pt}{0.800pt}}
\multiput(995.76,160.41)(-1.125,0.503){71}{\rule{1.985pt}{0.121pt}}
\multiput(999.88,157.34)(-82.881,39.000){2}{\rule{0.992pt}{0.800pt}}
\multiput(911.98,199.41)(-0.631,0.501){131}{\rule{1.209pt}{0.121pt}}
\multiput(914.49,196.34)(-84.491,69.000){2}{\rule{0.604pt}{0.800pt}}
\multiput(825.56,268.41)(-0.543,0.501){153}{\rule{1.070pt}{0.121pt}}
\multiput(827.78,265.34)(-84.779,80.000){2}{\rule{0.535pt}{0.800pt}}
\multiput(738.65,348.41)(-0.530,0.501){157}{\rule{1.049pt}{0.121pt}}
\multiput(740.82,345.34)(-84.823,82.000){2}{\rule{0.524pt}{0.800pt}}
\multiput(651.32,430.41)(-0.580,0.501){143}{\rule{1.128pt}{0.121pt}}
\multiput(653.66,427.34)(-84.659,75.000){2}{\rule{0.564pt}{0.800pt}}
\multiput(563.92,505.41)(-0.640,0.501){129}{\rule{1.224pt}{0.121pt}}
\multiput(566.46,502.34)(-84.461,68.000){2}{\rule{0.612pt}{0.800pt}}
\multiput(476.19,573.41)(-0.752,0.502){109}{\rule{1.400pt}{0.121pt}}
\multiput(479.09,570.34)(-84.094,58.000){2}{\rule{0.700pt}{0.800pt}}
\multiput(388.02,631.41)(-0.930,0.502){87}{\rule{1.681pt}{0.121pt}}
\multiput(391.51,628.34)(-83.511,47.000){2}{\rule{0.840pt}{0.800pt}}
\multiput(299.76,678.41)(-1.125,0.503){71}{\rule{1.985pt}{0.121pt}}
\multiput(303.88,675.34)(-82.881,39.000){2}{\rule{0.992pt}{0.800pt}}
\sbox{\plotpoint}{\rule[-0.200pt]{0.400pt}{0.400pt}}%
\put(221.0,123.0){\rule[-0.200pt]{293.416pt}{0.400pt}}
\put(1439.0,123.0){\rule[-0.200pt]{0.400pt}{157.549pt}}
\put(221.0,777.0){\rule[-0.200pt]{293.416pt}{0.400pt}}
\put(221.0,123.0){\rule[-0.200pt]{0.400pt}{157.549pt}}
\end{picture}
\caption{A plot of the error (defined as in \eqref{eq:experimentalErrorDefSoliton}) as a function of velocity. Note the exponential improvement in accuracy with velocity.
}
\end{figure}
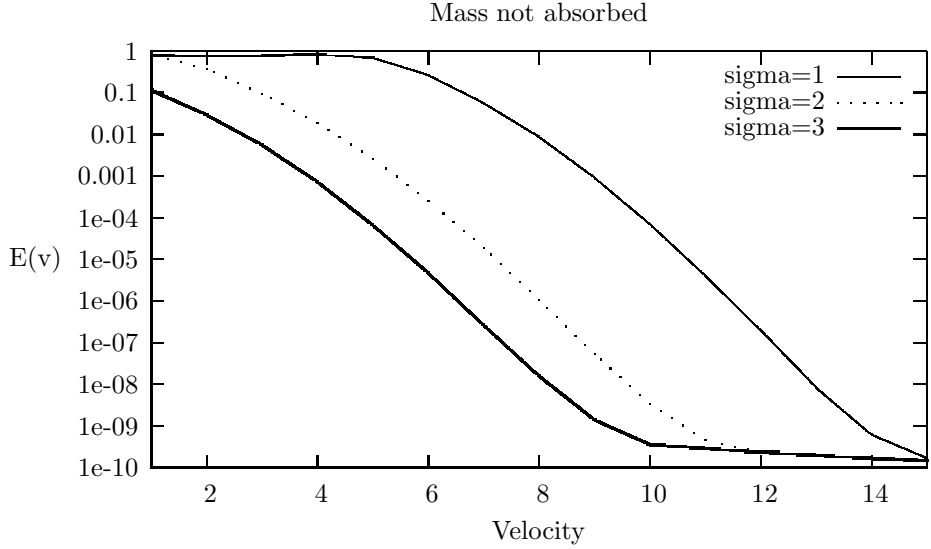

The TDPSF was used with $\Tstep=0.08/v$, $\xs = 0.20$, $\ks=2\pi/3.2$, and $\sdev=1,2,3$. We measured the following quantity:
\begin{equation}
  E(v)  =  \sup_{t < 200/v} \frac{
    \norm{\psi(x,t)-\psi_{ex}(x,t)}{L^2([-12,12])}
    }{
      \norm{\psi_{ex}(x,0)}{L^2(\mathbb{R})}
      } \label{eq:experimentalErrorDefSoliton}
\end{equation}
The function $\psi_{ex}(x,t)$ is the exact solution. The result of this experiment is plotted in Figure \ref{fig:solitonAbsorption}. The time $200/v$ was chosen since it is more than enough time for errors to return to the region $[-10,10]$. We believe the error floor near $10^{-10}$ visible in Figure \ref{fig:solitonAbsorption} is due to truncation errors in the calculation of the WFT, and errors due to time-stepping in solving the interior problem.

\begin{remark}
  \label{rem:compareToSzeftel}
  The paper \cite{szeftel:absorbingBoundaries} proposes an alternative method of absorbing boundaries (namely the paradifferential strategy), based on a novel method of approximating the Dirichlet-to-Neumann operator. A similar numerical test was performed for those boundary conditions. For a soliton at velocity $15$, Szeftel obtained $E(15)=0.08$ at best. For comparison, we obtain $E(15)=1.69\times 10^{-10}$ for $\sdev=1$ and $E(v)=1.40 \times 10^{-10}$ for $\sdev=3$. 

Our tests are not directly comparable. Our region of interest was $[-12,12]$ with TDPSF region on $[-25.6,-12]$ and $[12,25.6]$, as opposed to Szeftel who used $[-5,5]$ (although we used $256$ spatial lattice points rather than $200$, as used in \cite{szeftel:absorbingBoundaries}). We used FFT/Split Step propagation for the interior problem, as opposed to the finite differences of \cite{szeftel:absorbingBoundaries}. Additionally, storage of the history on the boundary was required in \cite{szeftel:absorbingBoundaries}, unlike the TDPSF.

It is quite surprising that the TDPSF beats the Dirichlet-to-Neumann boundary by such a large margin because \cite{szeftel:absorbingBoundaries}  takes the nonlinearity into account while the TDPSF assumes the nonlinearity is zero on the boundary.
\end{remark}

\section{Comparison to Other Methods}
\label{sec:comparisonToOtherMethods}

\subsection{Dirichlet-to-Neumann Boundaries}

The Dirichlet-to-Neumann map was originally constructed by Engquist and Majda in \cite{MR517938,MR0471386} (see also \cite{MR596431,MR658635}). Their guiding principle was that near the boundary, the physical optics approximation to wave flow is sufficiently accurate to filter off the outgoing waves.The TDPSF is a direct analogue of this - the Gaussian framelet elements behave (under the free flow) like classically free particles. We use a different method to filter, but the guiding principle is the same.

Modern approaches attempt to construct the exact solution on the boundary and then impose it as a boundary condition. In principle, this is the best possible approach, although in practice it is difficult to implement for dispersive equations. We briefly mention two major approaches that we are aware of, and remark that only one \cite{szeftel:absorbingBoundaries} even attempts to deal with nonlinear equations.

An additional problem is that this approach precludes the use of spectral methods (e.g. Algorithm \ref{algo:splitStep}) to solve the interior problem, reducing accuracy of the solution on the interior for dispersive equations. The FFT naturally imposes periodic boundaries rather than Dirichlet-to-Neumann, thus requiring the use of FDTD or other local methods.

\subsubsection{Current Constructions}

To deal with the free Schr\"odinger equation (no nonlinearity or potential), Lubich and Sch\"adle \cite{MR1924419,MR1785966} constructed an approximation to the exact integral kernel of the Dirichlet-Neumann operator based on a piecewise exponential (in time) approximation. This approach appears to work nicely for the free Schr\"odinger equation, although it is uncertain that it could be applied to the full Dirichlet-to-Neumann operator of a nonlinear equation or long range problem.

We are aware of only one fully nonlinear Dirichlet-to-Neumann operator, constructed by J. Szeftel in \cite{MR2114289}. Szeftel uses a modified version of the paradifferential calculus (see references in \cite{MR2114289}) in 1 space dimension to deal with smooth nonlinearities and potentials. He proves local well posedness of \eqref{eq:NLSE} of the boundary conditions assuming $H^6$ regularity of the initial data.

However, extensions to $\Rn$ appear highly nontrivial. The assumptions are significantly stronger than ours, and there are no error bounds. The numerical experiments look promising and the results appear accurate for radiative problems (see also Remark \ref{rem:compareToSzeftel}).

\subsection{Absorbing Potentials/ PML}

\subsubsection{Absorbing Potentials}

Absorbing (complex) potentials, described in \cite{neuhauser:complexPotentials}, are the current ``industry standard''. One can add a dissipative term $-i a(x) \psi(\vx,t)$ to the right side of \eqref{eq:NLSE} and solve it on the region $\FBox$. The function $a(x)$ is a positive function supported in $\FBox \setminus \IBox$. This has the effect of (partially) absorbing waves which have left the domain of interest, although it might create spurious reflections and dissipation. This approach is the mainstay of absorbing boundaries, being simple to program and compatible with spectral methods. 

The potential $a(x)$ must be tuned to the given problem. Given $\kmin$, $\kmax$, one must select the height and width of the absorber so that it kills most of the wave between $\kmin$ and $\kmax$. Waves with momentum lower than $\kmin$ are mostly reflected, and waves with momentum higher than $\kmax$ are mostly transmitted and wrap around the computational domain.

Heuristic calculations and numerical experiments suggest that the absorber must have width proportional to $C \kmax \ln(\epsilon) / \kmin$, with $C$ depending on the specific shape of the potential. The TDPSF works on a layer of width $C \ln(\epsilon) / \kmin$, which is smaller by a factor of $\kmax$. Spurious dissipation is an additional problem with absorbing potentials, as illustrated in Section \ref{sec:testsMediumRange}.

\subsubsection{Perfectly Matched Layers}

Perfectly Matched Layers (PML) are a variation on this approach, first proposed in \cite{MR1294924} for Maxwell's equations. In \cite{MR2032866} they are designed and tested for the $1$ dimensional free Schr\"odinger equation, with reasonable results.  

In the special case of the Schr\"odinger equation\footnote{In general, one constructs the PML by exterior complex scaling, but for the Schr\"odinger equation it takes this simple form.}, the PML consists of adding a term $-i a(x) \Lap \psi(\vx,t)$ to the right side of \eqref{eq:NLSE} instead of merely $-i a(\vx) \psi(\vx,t)$. If $a(\vx)$ is chosen carefully, one can completely eliminate reflections at the interface (the boundary of $\operatorname{supp} a(\vx)$) for certain frequencies. It also has the property of dissipating high frequency waves more strongly than low frequency ones, thereby requiring a boundary layer of size only $C \ln(\epsilon)/\kmin$. 

As discussed in remark \ref{remark:commentsOnPML}, numerical tests using the PML of \cite{MR2032866} achieved errors of only $10^{-3}$. The PML is likely to have the same problem as complex absorbing potentials with spurious dissipation. Some PML methods are unstable \cite{abarbanel:PMLinstability}. The PML method for the Schr\"odinger equation is still very much undeveloped, so a more detailed comparison is difficult to make.

\section{Conclusion}
\label{sec:outlook}

We have described in this work a new method, the Time Dependent Phase Space Filter, of filtering outgoing waves for the nonlinear Schr\"odinger equation. Unlike absorbing potentials and (most likely) PML methods, the TDPSF filters only those regions of phase space containing outgoing waves. Dissipative terms on the boundary filter all waves, including waves which should have returned to the computational region from the boundary (see Section \ref{sec:testsMediumRange}).

Our method is easier to construct than the Dirichlet-Neumann map, being based on standard techniques of signal processing rather than complicated pseudo/para-differential calculus. It is local in time, whereas the Dirichlet-Neumann map is history dependent. The TDPSF is also compatible with spectral methods, unlike the Dirichlet-to-Neumann map. 

Unlike all other methods we are aware of, the TDPSF fails gracefully, in the sense that when it cannot filter the outgoing waves it notifies the user rather than reporting an incorrect calculation as correct.

The main problem with the method is the problem of filtering low frequency outgoing waves, a problem shared with all other methods. We are currently working on an extension of the TDPSF \cite{us:multiscale} which will reduce the cost of filtering from $C \ln(\epsilon )/\kmin$ to $C \ln (\epsilon) \ln (\kmin)$. 

{\bf Acknowledgements:} We thank O. Costin and M. Kiessling for useful discussions, and R. Falk for reading the manuscript. This work was supported by NSF grant DMS01-00490.

\nocite{MR818063}
\nocite{MR898052}
\nocite{MR966733}
\nocite{MR0436612}
\nocite{MR0471386}
\nocite{MR517938}

\nocite{MR1819643}
\nocite{MR2032866}

\bibliographystyle{plain}
\bibliography{../stucchio.bib}

\end{document}